\documentclass[11pt,lefteqn]{article}

\usepackage{amsmath,amsfonts,amsthm,amssymb,graphics,epsfig,color}

\newtheorem{theorem}{Theorem}

\newtheorem{corollary}[theorem]{Corollary}

\begin{document}
\date {\today}
\title{When is the Resolvent Like a Rank One Matrix?}

\author{Anne Greenbaum\footnote{University of Washington, Applied Math Dept., Box 353925, Seattle, WA 98195.  
email:  greenbau@uw.edu}, Faranges Kyanfar\footnote{Shahid Bahonar University of Kerman, Applied Math Dept., 
Kerman, Iran.  kyanfar@uk.ac.ir}, Abbas Salemi\footnote{Mahani Math Center, 
Afzalipour Research Institute, Shahid Bahonar University of Kerman,
Kerman, Iran.  salemi@uk.ac.ir}}

\maketitle

\begin{abstract}
For a square matrix $A$, the resolvent of $A$ at a point $z \in \mathbb{C}$
is defined as $(A-zI )^{-1}$.  We consider the set of points $z \in \mathbb{C}$
where the relative difference in 2-norm between the resolvent and the nearest
rank one matrix is less than a given number $\epsilon \in (0,1)$.  We establish a 
relationship between this set and the $\epsilon$-pseudospectrum of $A$,
and we derive specific results about this set for Jordan blocks and for
a class of large Toeplitz matrices.  
We also derive disks about the eigenvalues of $A$ that are contained in this set,  
and this leads to some new results on disks about the eigenvalues that are
contained in the $\epsilon$-pseudospectrum of $A$.  In addition, we consider the set of
points $z \in \mathbb{C}$ where the absolute value of the inner product of
the left and right singular vectors corresponding to the largest singular
value of the resolvent is less than $\epsilon$.  We demonstrate numerically
that this set can be almost as large as the one where the relative difference
between the resolvent and the nearest rank one matrix is less than $\epsilon$
and we give a partial explanation for this.  Some possible applications are discussed.
\end{abstract}

\section{Introduction} \label{sec:intro}
Let $A$ be an $n$ by $n$ matrix, $n \geq 2$.  The {\em resolvent} of $A$ at a point $z \in \mathbb{C}$
is defined as $(A-zI )^{-1}$.  If the singular values of $(A - zI )^{-1}$
are denoted as $\sigma_1 (z) \geq \cdots \geq \sigma_n (z)$, and the corresponding
left and right singular vectors are $u_j (z)$ and $v_j (z)$, $j=1, \ldots , n$,
then the nearest rank one matrix to $(A-zI )^{-1}$ in the 2-norm is 
$\sigma_1 (z) u_1 (z) v_1 (z )^H$, and the 2-norm of the difference between
the resolvent and this rank one matrix is $\sigma_2 (z)$.  The {\em relative difference}
in 2-norm between the resolvent and this rank one matrix can be defined as
\begin{equation}
\| (A-zI )^{-1} - \sigma_1 (z) u_1 (z) v_1 (z )^H \|_2 / \| (A - zI )^{-1} \|_2 =
\sigma_2 (z) / \sigma_1 (z) . \label{reldiff}
\end{equation}
We are interested in the set $S_{\epsilon} (A)$ where this relative difference
is less than a given $\epsilon \in (0,1)$:
\begin{equation}
S_{\epsilon} (A) := \{ z \in \mathbb{C} : \sigma_2 (z) / \sigma_1 (z) < \epsilon \} .
\label{Seps}
\end{equation}
If $z$ is an eigenvalue of $A$, we will define the ratio in (\ref{Seps}) to be $0$,
so the spectrum of $A$ is always included in $S_{\epsilon} (A)$.
If $A$ is a multiple of the identity, then this is all that is included in 
$S_{\epsilon} (A)$.

The following theorem lists some simple properties of $S_{\epsilon} (A)$.
\begin{theorem}
\label{thm:simple}
Let $A$ be a square matrix and let $S_{\epsilon} (A)$ be defined by (\ref{Seps}) for 
$\epsilon \in (0,1)$.
Then:
\begin{enumerate}
\item
For any constant $c$, $S_{\epsilon} (A + cI ) = c + S_{\epsilon} (A)$.
\item
For any constant $c$, $S_{\epsilon} (cA) = c S_{\epsilon} (A)$.
\item
$S_{\epsilon} ( A^H ) = \overline{S_{\epsilon} (A)}$.
\item
For any unitary matrix $Q$, $S_{\epsilon} ( Q^H A Q ) = S_{\epsilon} (A)$.
\item
$S_{\epsilon} (A)$ is the set of points $z \in \mathbb{C}$ for which the ratio of the smallest
to the second smallest singular value of $A-zI$ is less than $\epsilon$.  
(If $A-zI$ has a zero singular value, we define this ratio to be zero, even if $A-zI$
has more than one zero singular value.)
\item
If $A$ is {\em block diagonal}, $A := \mbox{blockdiag}( A_1 , \ldots , A_m )$, then 
$S_{\epsilon} (A) \subseteq \bigcup_{j=1}^m S_{\epsilon} ( A_j )$.
In the other direction $S_{\epsilon} (A)$ contains the set of points $z$ such that
$z \in S_{\epsilon} ( A_k )$ for some block $k$, and for each $\ell \neq k$,
the smallest singular value of $A_{\ell} - zI$ is greater than or equal to the
second smallest singular value of $A_k - zI$. 
\end{enumerate}
\end{theorem}
\begin{proof}
\begin{enumerate}
\item
Since $(A-zI )^{-1} = (A+cI - (z+c)I )^{-1}$, $z \in S_{\epsilon} (A)$ if and only if
$z+c \in S_{\epsilon} (A+cI)$.
\item
Since the ratio in (\ref{Seps}) remains unchanged if $A$ is replaced by $cA$ and $z$ is
replaced by $cz$, it follows that $S_{\epsilon} (cA) = c S_{\epsilon} (A)$.
\item
Since $( A^H - \bar{z} I )^{-1} = ( (A-zI )^{-1} )^H$, the singular values
of $( A^H - \bar{z} I )^{-1}$ are the same as those of $(A-zI )^{-1}$, and therefore
$\bar{z} \in S_{\epsilon} ( A^H )$ if and only if $z \in S_{\epsilon} (A)$.
\item
Since the singular values of $( Q^H A Q - zI )^{-1} = Q^H (A-zI )^{-1} Q$ are the
same as those of $(A-zI )^{-1}$, it follows that $S_{\epsilon} ( Q^H A Q ) = S_{\epsilon} (A)$.
\item
The singular values of $A-zI$ are the inverses of the singular values of $(A - zI )^{-1}$,
so the ratio of the smallest to the second smallest singular value of $A-zI$ is the same
as the ratio of the second largest to the largest singular value of $(A-zI )^{-1}$.
\item
Using result 4, the singular values of $A-zI$ are the singular values of all of the blocks
$A_j - zI$, $j=1, \ldots , m$.  Suppose the smallest singular value of $A-zI$ is the 
smallest singular value of $A_k - zI$.  The second smallest singular value of $A-zI$ is
less than or equal to the second smallest singular value of $A_k - zI$, so that if
$z \in S_{\epsilon} (A)$, then $z \in S_{\epsilon} ( A_k )$.
In the other direction, if the singular values of all blocks $A_{\ell} - zI$, $\ell \neq k$
are all greater than or equal to the two smallest singular values of $A_k - zI$, 
then the ratio of the smallest to the second smallest singular value of $A-zI$
is the same as that of $A_k - zI$, so if $z \in S_{\epsilon} ( A_k )$, then
$z \in S_{\epsilon} (A)$. 
\end{enumerate}
\end{proof}

In Section \ref{sec:Numerical}, we present numerical examples showing that 
when $A$ has ill-conditioned
eigenvalues, the set $S_{\epsilon} (A)$ can be large -- much larger than, say,
a union of disks of radius $\epsilon$ about the eigenvalues of $A$.

This set clearly is related to the 2-norm {\em $\epsilon$-pseudospectrum} 
of $A$ \cite{TrefEmb}, which can be defined as 
\begin{equation}
\Lambda_{\epsilon} (A) := \{ z \in \mathbb{C} : 1 / \sigma_1 (z) < \epsilon \} .
\label{epspseudo}
\end{equation}
Note a difference in scaling, however, in that, unlike property 2 in Theorem \ref{thm:simple},
it is not $\Lambda_{\epsilon} (cA)$ that is equal to $c \Lambda_{\epsilon} (A)$ 
but $\Lambda_{|c| \epsilon} (cA) = c \Lambda_{\epsilon} (A)$.
The relationship between the two sets is characterized in Section \ref{sec:pseudo}.

Sections \ref{sec:Jordan} and \ref{sec:Toeplitz} deal with Jordan blocks
and large Toeplitz matrices.  For a Jordan block, we derive a disk about
the eigenvalue that is contained in $S_{\epsilon} (A)$, while for more general
Toeplitz matrices, we discuss what happens when the {\em symbol} of the 
Toeplitz operator has winding number $\pm 1$ about the origin.  For any shift
$-z$ that maintains this property, we use the {\em splitting property} to
argue that as $n \rightarrow \infty$,
the smallest singular value of $A-zI$ approaches $0$ and the second smallest
singular value approaches a nonzero value.  Thus, for $n$ sufficiently large,
$z$ will lie in $S_{\epsilon} (A)$.

Section \ref{sec:genl} is devoted to more general diagonalizable matrices,
and we again derive disks about the eigenvalues that are contained in $S_{\epsilon} (A)$.
This leads to a new result about disks about the eigenvalues that are 
contained in the $\epsilon$-pseudospectrum.

It was observed in \cite{GreenWell} that not only can $S_{\epsilon} (A)$
be a large set, but that throughout much of this set it may be the case that
$| u_1 (z )^H v_1 (z) | < \epsilon$.  If $z$ is extremely close to an
ill-conditioned eigenvalue $\lambda$ of $A$ ({\em much} closer to $\lambda$ 
than to any other eigenvalue of $A$), then this is to be expected since
$u_1 (z)$ and $v_1 (z)$ will be close to the unit left and right eigenvectors 
$x$ and $y$ associated with the eigenvalue $\lambda$; the {\em condition number}
of $\lambda$ is defined as $1/| y^H x |$, and $\lambda$ is said to be ill-conditioned
when $| y^H x |$ is tiny.  It was observed, however, that 
$u_1 (z)$ and $v_1 (z)$ may be nearly orthogonal to each other even when $z$
is too far away from any eigenvalue of $A$ for $u_1 (z)$ and $v_1 (z)$ to 
resemble left and right eigenvectors.  
In Section \ref{sec:diffeq}, we give a partial explanation for this observation.

Finally, in Section \ref{sec:open} we mention some possible applications and remaining
open problems.

Throughout the paper, superscript $^H$ denotes the complex conjugate transpose
of a vector or matrix and $\| \cdot \|$ denotes the 2-norm for vectors and the
corresponding spectral norm (i.e., the largest singular value) for matrices.
In this section, we have denoted the singular values of $(A-zI )^{-1}$ as
$\sigma_j (z)$, $j=1, \ldots , n$, to emphasize the fact that $A$ is fixed and
$\sigma_j$ varies with $z$, and we will retain this notation in Section \ref{sec:Numerical} 
and in Section \ref{sec:open}.
In Section \ref{sec:diffeq}, we deal with matrices of the form $A-r e^{i \theta} I$, where $A$
and $\theta$ are fixed, and we denote their singular values as $\sigma_j (r)$.  
Throughout most of the paper, however, we deal with a number of different matrices 
so we will denote the singular values of a matrix $B$ explicitly as $\sigma_j (B)$.

\section{Numerical Examples} \label{sec:Numerical}
\begin{figure}[ht]
\centerline{\epsfig{file=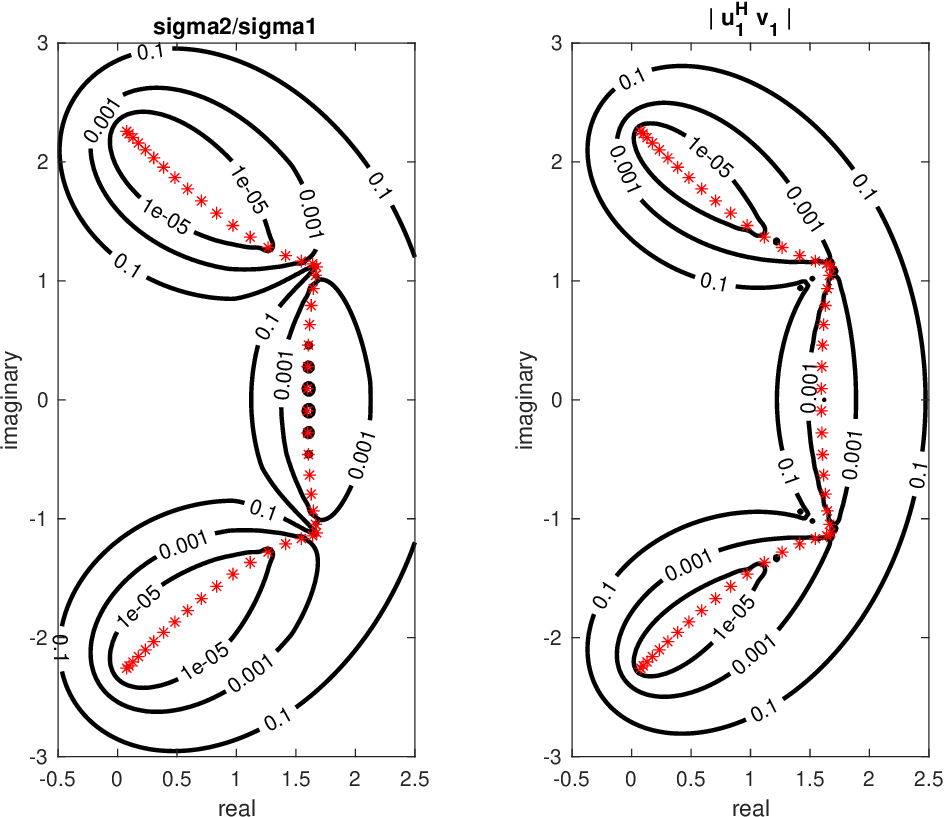,width=2.5in}}
\caption{Contour plots of ratios of second largest to largest singular value of the resolvent (left) and of
absolute value of inner product of left and right singular vectors corresponding to largest singular value
of the resolvent (right).  Matrix is the Grcar matrix of size $n=50$.
Eigenvalues are shown with red asterisks.}
\label{fig:grcar50}
\end{figure}

\begin{figure}[ht]
\centerline{\epsfig{file=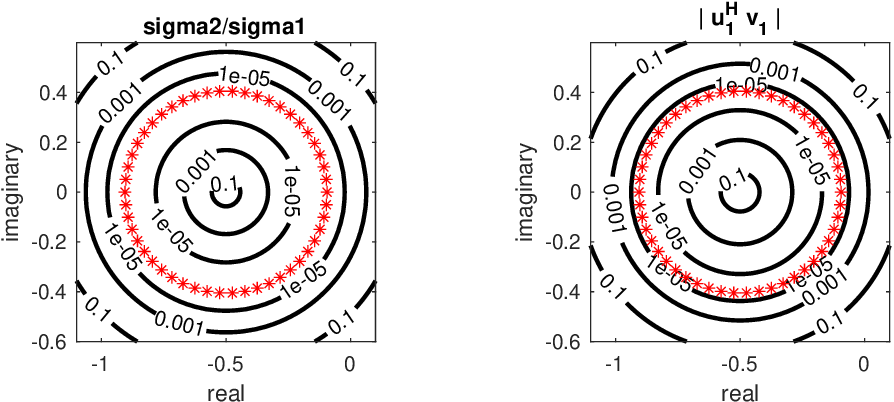,width=4in}}
\caption{Contour plots of ratios of second largest to largest singular value of the resolvent (left) and of
absolute value of inner product of left and right singular vectors corresponding to largest singular value
of the resolvent (right).  Matrix is the transient\_demo matrix of size $n=50$.  
Eigenvalues are shown with red asterisks.}
\label{fig:transient50}
\end{figure}

\begin{figure}[ht]
\centerline{\epsfig{file=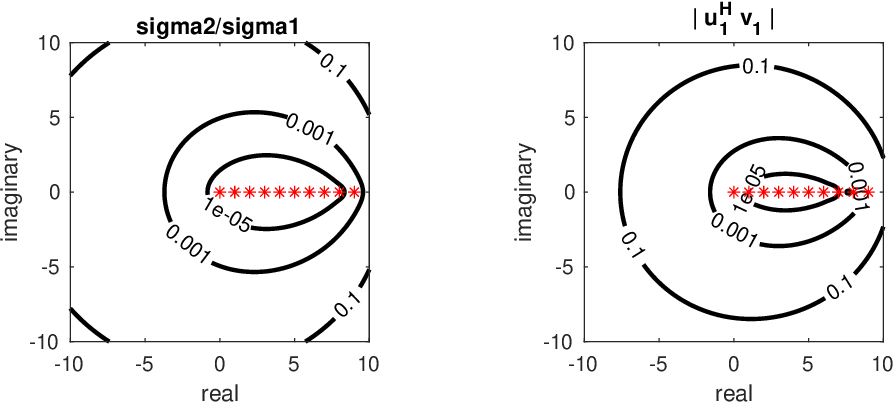,width=4in}}
\caption{Contour plots of ratios of second largest to largest singular value of the resolvent (left) and of
absolute value of inner product of left and right singular vectors corresponding to largest singular value
of the resolvent (right).  Matrix is the sampling matrix of size $n=10$.  
Eigenvalues are shown with red asterisks.}
\label{fig:sampling10}
\end{figure}

The phenomenon described in the Introduction is illustrated in Figures 
\ref{fig:grcar50}--\ref{fig:sampling10} 
for three matrices with ill-conditioned eigenvalues --
the Grcar matrix\footnote{In MATLAB, gallery('grcar',50).  
This matrix has $-1$'s on the subdiagonal, $1$'s on the main diagonal 
and the first three super-diagonals, and $0$'s elsewhere.}
of size $n=50$, the \verb+transient_demo+ matrix\footnote{In \cite{eigtool},
transient\_demo(50).} of size $n=50$, and the sampling matrix\footnote{In
MATLAB, gallery('sampling',10).  This matrix has integer eigenvalues,
$0,1, \ldots , 9$, that are ill-conditioned.} of size $n=10$.  
For the Grcar matrix, the eigenvalue condition numbers vary
from about $1.5e+2$ to $2.2e+7$, with the most ill-conditioned eigenvalues 
being those with imaginary part near $\pm 2$ and the best conditioned being 
those with real part about $1.6$ and imaginary part $\pm 1.1$ (those near 
where the curved sections of eigenvalues meet).  All eigenvalues of
the transient\_demo matrix have the same condition number, about $7.2e+5$.
The eigenvalues of the sampling matrix have condition numbers ranging
from about $1.4e+3$ to $1.3e+6$, with the interior eigenvalues being
more ill-conditioned than those at the ends.  The plots on the
left show contours of the ratio of the second largest to the
largest singular value of $(A - zI )^{-1}$, which, as noted previously,
is the relative difference in 2-norm between $(A-zI )^{-1}$
and the nearest rank one matrix.  

\medskip

\medskip

\medskip

The right plots in the figures show the
absolute values of the inner products of the left and right singular vectors
$u_1 (z)$ and $v_1 (z)$ corresponding to the largest singular value 
of $(A-zI )^{-1}$.  It can be seen from the figures that these 
inner products are small even when $z$ is too far away from any eigenvalue 
to argue that $u_1 (z)$ and $v_1 (z)$ resemble unit left and right eigenvectors. 

\begin{figure}[ht]
\centerline{\epsfig{file=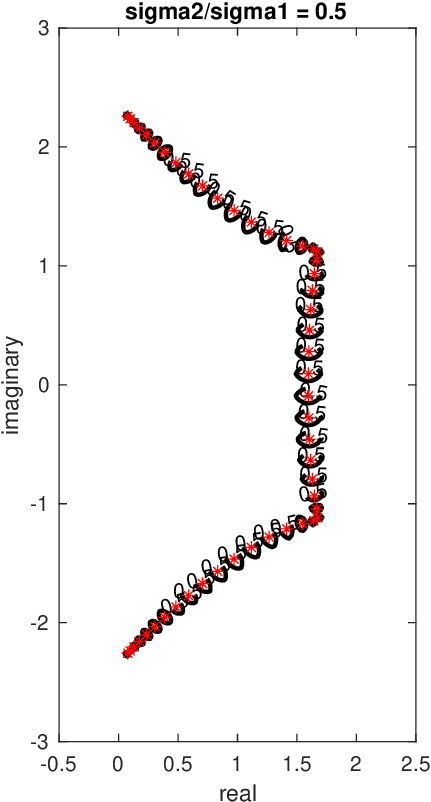,width=1.2in}}
\caption{Contour plot of ratios of second largest to largest singular value 
of the resolvent for a normal matrix with the same eigenvalues
as the Grcar matrix of size $n=50$.  The $0.5$ contour is plotted
because contours with much smaller values are barely visible on the
scale of the graph.
Eigenvalues are shown with red asterisks.}
\label{fig:normalgrcar50}
\end{figure}

In contrast, if we take a normal matrix whose eigenvalues are the
same as those of, say, the Grcar matrix, then the contours plotted
in Figure \ref{fig:grcar50} will not be visible on the scale of the graph; 
for illustration, we plot in Figure \ref{fig:normalgrcar50} the contour where 
$\sigma_2 (z) / \sigma_1 (z) = 0.5$
for a diagonal matrix whose eigenvalues are those
of the Grcar matrix.  Clearly, the phenomenon illustrated in 
Figures \ref{fig:grcar50}--\ref{fig:sampling10} is limited to
matrices with at least some ill-conditioned eigenvalues.

\section{Relation to Pseudospectra} \label{sec:pseudo}
\begin{figure}[ht]
\centerline{\epsfig{file=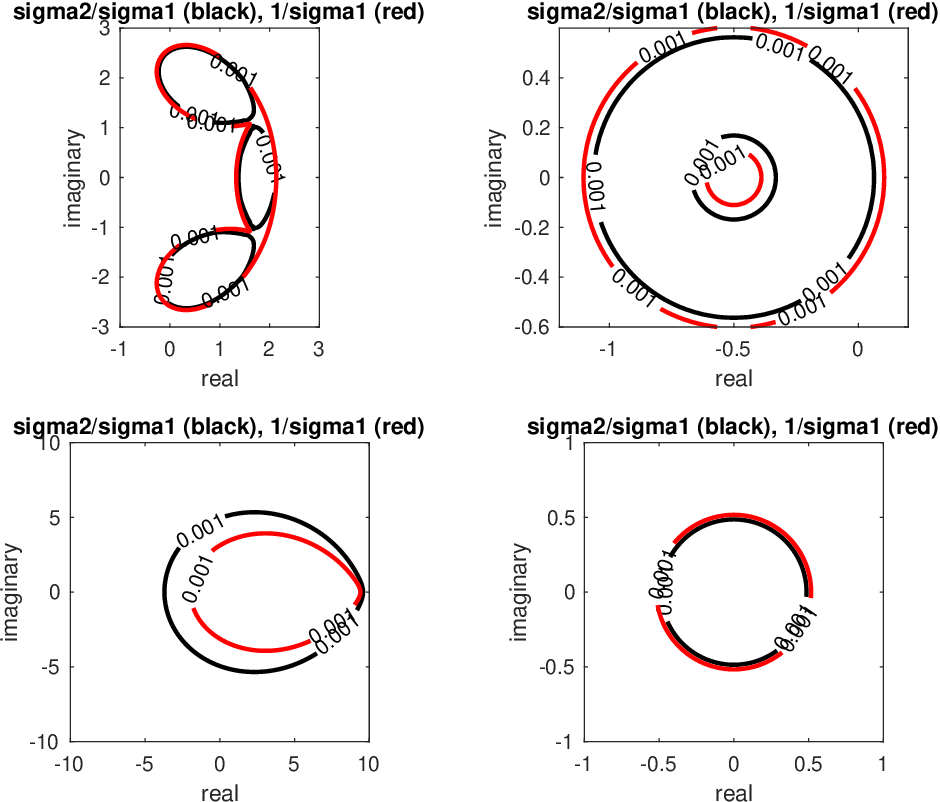,width=3in}}
\caption{Comparison of curve where $\sigma_2 ( (A-zI )^{-1}) / \sigma_1 ( (A-zI )^{-1}) =
10^{-3}$ (black) to that where $1 / \sigma_1 ( (A-zI )^{-1}) = 10^{-3}$ (red); i.e., to
the boundary of the $10^{-3}$-pseudospectrum.  Upper left is the Grcar matrix
of size $n=50$, upper right is the transient\_demo matrix of size
$n=50$.  Lower left is the sampling matrix of size $n=10$, lower right
is a Jordan block of size $n=10$.}
\label{fig:pseudo}
\end{figure}

The $\epsilon$-pseudospectrum \cite{TrefEmb} of a matrix $A$ is usually defined as
\begin{equation}
\{ z \in \mathbb{C} : \| (A-zI )^{-1} \| > \epsilon^{-1} \} .
\label{pseudo}
\end{equation}
For the 2-norm this is clearly equivalent to (\ref{epspseudo}), so
the boundary of the $\epsilon$-pseudospectrum is the
curve on which $1/ \sigma_1 ( (A-zI )^{-1} ) = \epsilon$, while the
left plots of Figures \ref{fig:grcar50}--\ref{fig:sampling10} show
the curves on which $\sigma_2 ( (A-zI )^{-1} ) / \sigma_1 ( (A-zI )^{-1} ) =
\epsilon$, for $\epsilon = 10^{-5},~10^{-3},~10^{-1}$.  In Figure
\ref{fig:pseudo}, we compare the $10^{-3}$-pseudospectrum with
the region in which $\sigma_2 ( (A-zI )^{-1} ) / \sigma_1 ( (A-zI )^{-1} ) < 10^{-3}$
for each of the three matrices in Figures 
\ref{fig:grcar50}--\ref{fig:sampling10}, as well as for a
$10$ by $10$ Jordan block.  For these matrices and for $z$ in the regions plotted, 
$\sigma_2 ( (A-zI )^{-1} )$ is not too different from $1$, so the curves look similar.  

Recall from (\ref{Seps}) and (\ref{epspseudo}), however, that there is a difference 
in scaling between $S_{\epsilon} (A)$ and $\Lambda_{\epsilon} (A)$:
For a constant $c$, $S_{\epsilon} (cA) = c S_{\epsilon} (A)$ but 
$\Lambda_{|c| \epsilon} (cA) = c \Lambda_{\epsilon} (A)$.
To compare $S_{\epsilon} (A)$ to $\Lambda_{\epsilon} (A)$, we should multiply
$A$ by a constant so that, say, the mean value of $\sigma_2 (( A-zI )^{-1} )$ on
$\Lambda_{\epsilon} (A)$ is around $1$.

To quantify the relation with pseudospectra, we will work with the matrix
$A-zI$, whose singular values are the inverses of those of $(A-zI )^{-1}$
and whose right/left singular vectors are the same as the left/right singular
vectors of $(A - zI )^{-1}$.  As noted in item 5 of Theorem \ref{thm:simple},
the ratio of the smallest to the second smallest
singular value of $A-zI$ is the same as the ratio of the second largest to the
largest singular value of $(A-zI )^{-1}$. 
We will use the following inequality, which is a consequence of the Wielandt-Hoffman 
theorem and can be found, for instance, in [Corollary 8.6.2]\cite{GolubvanLoan}:  
The decreasingly ordered singular values of two $n$ by $n$ matrices $B$ and $B+E$ 
satisfy
\begin{equation}
| \sigma_j (B+E) - \sigma_j (B) | \leq \sigma_1 (E) . \label{WieHoff}
\end{equation}
With $E=-zI$ and $j=n-1$, this implies that
\begin{equation}
\sigma_{n-1} (B) + |z| \geq \sigma_{n-1} (B-zI) \geq \sigma_{n-1} (B) - |z| . 
\label{sigmanm1bound}
\end{equation}
This leads to the following theorem relating the $\epsilon_1$-pseudospectrum of $A$
to the region where $\sigma_2 ( (A-zI )^{-1} ) / \sigma_1 ( (A-zI )^{-1} ) < \epsilon_2$,
for what may be nearby values $\epsilon_1$ and $\epsilon_2$.

\begin{theorem}
\label{thm:relation}
Let $A$ be an $n$ by $n$ matrix with distinct eigenvalues $\lambda_1 , \ldots , \lambda_n$,
and let $\Lambda_{\epsilon} (A)$ be defined by (\ref{epspseudo}) and $S_{\epsilon} (A)$
be defined by (\ref{Seps}).
Consider the matrices $A_{0,j} := A - \lambda_j I$, $j=1, \ldots , n$, each of which 
has a zero singular value, $\sigma_n ( A_{0,j} )$.  Let $\sigma_{n-1} ( A_{0,j} ) > 0$ 
denote the second smallest singular value of $A_{0,j}$.
For $r > 1$, $S_{\epsilon (r/( (r-1) \min_j \sigma_{n-1} ( A_{0,j} )))} (A)$ contains 
\begin{equation}
\Lambda_{\epsilon} (A) \cap \left( \bigcup_{j=1}^n \mathbb{D} ( \lambda_j , 
\sigma_{n-1} ( A_{0,j} )/r ) \right) , \label{pseudointersection}
\end{equation}
where $\mathbb{D}(c,R)$ denotes the open disk about $c$ of radius $R$.
In the other direction, 
\begin{equation}
S_{\epsilon} (A) \cap
\left( \bigcup_{j=1}^n \mathbb{D} ( \lambda_j , \sigma_{n-1} ( A_{0,j} )/r ) \right)
\label{newsetintersection}
\end{equation}
is contained in $\Lambda_{\epsilon ( ( \max_{j} \sigma_{n-1} ( A_{0,j} ) ) ((r+1)/r))} 
(A)$.
\end{theorem}
\begin{proof}
For $z \in \Lambda_{\epsilon} (A)$, $\sigma_n (A-zI) < \epsilon$.  
From (\ref{sigmanm1bound}), the second smallest singular value of 
$A - zI = A_{0,j} - (z- \lambda_j )I$ is greater than or equal to 
$\sigma_{n-1} ( A_{0,j} ) - |z - \lambda_j |$.  Hence for $z \in 
\mathbb{D} ( \lambda_j, \sigma_{n-1} ( A_{0,j} )/r )$, the second smallest singular value 
of $A - zI$ is greater than $\sigma_{n-1} ( A_{0,j} ) ( r-1 )/r$.  
It follows that for $z$ in $\Lambda_{\epsilon} (A) \cap 
\mathbb{D} ( \lambda_j , \sigma_{n-1} ( A_{0,j} ) / r )$, we have
$\sigma_n (A-zI) / \sigma_{n-1} (A-zI) < \epsilon (r / ((r-1) \sigma_{n-1} ( A_{0,j} ) ))$.
Taking the max over $j$ of the bounds on $\sigma_n (A-zI) / \sigma_{n-1} (A-zI)$
establishes (\ref{pseudointersection}).  

Now suppose $\sigma_n (A-zI) / \sigma_{n-1} (A-zI) < \epsilon$ and $z$ lies 
in one of the disks $\mathbb{D} ( \lambda_j , \sigma_{n-1} ( A_{0,j}  ) / r )$.
Then we can write
\[
\sigma_n (A-zI) < \epsilon \sigma_{n-1} (A-zI) = \epsilon \sigma_{n-1}
( A_{0,j} - (z- \lambda_j )I ) \leq \epsilon ( \sigma_{n-1} ( A_{0,j} ) + 
|z- \lambda_j | ) < \epsilon \sigma_{n-1} ( A_{0,j} ) (r+1)/r .
\]
Taking the max over $j$ of the bounds on $\sigma_n (A-zI)$ 
establishes (\ref{newsetintersection}).
\end{proof}

{\em If} the second smallest singular value of each matrix $A - \lambda_j I$ is
on the order of $1$, as it is for the matrices used in Figure \ref{fig:pseudo},
and if $\epsilon$ is small enough so that $\Lambda_{\epsilon} (A)$ and
$S_{\epsilon} (A)$ lie mostly inside the union of disks in (\ref{pseudointersection})
and (\ref{newsetintersection}), then Theorem \ref{thm:relation} establishes a close
relationship between pseudospectra and sets where the resolvent is close
to a rank one matrix.  Taking $r=2$, for example, if $\min_j \sigma_{n-1} ( A_{0,j} ) =
\max_j \sigma_{n-1} ( A_{0,j} ) = 1$ and assuming that $\Lambda_{\epsilon} (A)$ 
and $S_{\epsilon} (A)$ lie inside the union of disks in (\ref{pseudointersection})
and (\ref{newsetintersection}), the theorem shows that 
$\Lambda_{\epsilon /2} (A) \subset S_{\epsilon} (A) \subset \Lambda_{(3/2) \epsilon} (A)$
and $S_{(2/3) \epsilon} (A) \subset \Lambda_{\epsilon} (A) \subset S_{2 \epsilon} (A)$.

If the second smallest singular value of each matrix $A_{0,j}$, $j=1, \ldots , n$,
is {\em not} on the order of $1$, then this close relationship may not hold.
For example, suppose one starts with the transient\_demo matrix (whose eigenvalues
all have the same condition number), computes an eigendecomposition
$X \Lambda X^{-1}$, then randomly permutes the diagonal entries of $\Lambda$
and forms the matrix $A = X \mbox{diag}( \lambda_{\pi_1} , \ldots , \lambda_{\pi_n} ) 
X^{-1}$.  Then $A$ will have the same eigenvalues with the same condition numbers and 
the same eigenvectors as the transient\_demo matrix, the only difference being that 
different eigenvalues are associated with different eigenvectors.  However, the 
matrices $A_{0,j} = A - \lambda_j I$ in Theorem \ref{thm:relation}, are entirely different.  
In our experiments with $n=50$, the matrices $A_{0,j}$ had several singular values
in the $10^{-8}$ range.  In this case,
the sets $\Lambda_{\epsilon} (A)$ and $S_{\epsilon} (A)$ are entirely different.

\section{Analysis for a Jordan Block} \label{sec:Jordan}
To understand more about why Figures \ref{fig:grcar50}--\ref{fig:sampling10} 
show such large regions of the complex plane where the resolvent,
$(A-zI )^{-1}$, closely resembles a rank one matrix 
we will look for disks about the eigenvalues of $A$ that can
be shown to be contained in these regions.  Since 
$\sigma_n (A-zI) = \min_{w \neq 0} \| (A-zI) w \| / \| w \|$, one can obtain
a good upper bound on $\sigma_n (A-zI)$ by displaying a vector $w$ that approximately
minimizes this ratio.  Here we do this for an $n$ by $n$ Jordan block.

\begin{theorem}
\label{thm:Jordan}
Let $A$ be an $n$ by $n$ Jordan block so that
\[
A - zI = \left[ \begin{array}{cccc}
-z & 1      &        &   \\
   & \ddots & \ddots &   \\
   &        & \ddots & 1 \\
   &        &        & -z \end{array} \right]_{n \times n} .
\]
The smallest singular value $\sigma_n (A-zI)$ is less than or equal to $| z |^n$,
and the second smallest singular value $\sigma_{n-1} (A-zI)$ is greater than or 
equal to $1 - |z|$. It follows that for any $\epsilon \in (0,1)$,
the set of points $z$ at which $\sigma_n (A-zI) / \sigma_{n-1} (A-zI) <
\epsilon$ contains a disk about the origin of radius $r$,
where $r^n / (1-r) = \epsilon$; this always includes a disk of radius
$\frac{1}{2} \epsilon^{1/n}$.
\end{theorem}

\begin{proof}
If $w = [1, z, z^2 , \ldots , z^{n-1} ]^T$, then $(A-zI)w =
[0, \ldots , 0, - z^n ]^T$.  Since $\sigma_n (A-zI) =
\min_{w \neq 0} \| (A-zI) w \| / \| w \|$, it follows that
\begin{equation}
\sigma_n (A-zI ) \leq \frac{| z^n |}{\sqrt{\sum_{j=0}^{n-1} | z^j |^2}} \leq 
| z^n | . \label{sigmanbound}
\end{equation}

The Jordan block $A$ has a zero singular value and $n-1$ ones as singular
values, since $A^H A = \mbox{diag} ( 0, 1, \ldots , 1 )$.  
Combining (\ref{sigmanbound}) and (\ref{sigmanm1bound}), we see that
\begin{equation}
\sigma_n (A-zI) / \sigma_{n-1} (A-zI) \leq | z |^n / (1- |z| ) . \label{boundn}
\end{equation}
It follows that $z \in S_{\epsilon} (A)$ if $| z |^n / (1- |z|) < \epsilon$.

Finally, given $\epsilon \in (0,1)$, if $|z| < \frac{1}{2} \epsilon^{1/n}$,
then $|z |^n < 2^{-n} \epsilon$ and $1- |z| > 1/2$ so that
\[
\sigma_n (A-zI) / \sigma_{n-1} (A-zI) < 2^{-(n-1)} \epsilon \leq 
\epsilon .
\]
\end{proof}

Note that the radius of the black circle in the lower right plot of Figure \ref{fig:pseudo}
is about $0.486$, and the $10$th root of $10^{-3}$ is about $0.501$.  Theorem \ref{thm:Jordan}
guarantees that $S_{10^{-3}} (A)$ contains a disk about the origin of radius $r$, where
$r^{10} / (1-r) = 10^{-3}$, and solving for $r$ we find $r \approx 0.470$.

This can be compared with a known result about the $\epsilon$-pseudospectrum
of a Jordan block.  It is shown in \cite{ReichTref} that the $\epsilon$-pseudospectrum
of an $n$ by $n$ Jordan block contains a disk about the origin of radius $\epsilon^{1/n}$ 
and is contained in a disk about the origin of radius $1+ \epsilon$.  The lower bound 
is demonstrated to be a much sharper estimate of the actual radius.

\section{Large Toeplitz Matrices} \label{sec:Toeplitz}
For more general large Toeplitz matrices, results about how close the resolvent is
to a rank one (or otherwise low rank) matrix can be derived from what is known
as the {\em splitting phenomenon} \cite[Theorem 9.4]{BottGrud}.
An infinite Toeplitz matrix has the form
\[
\left[ \begin{array}{cccc} a_0 & a_{-1} & a_{-2} & \ldots \\
a_1 & a_0 & a_{-1} & \ldots \\ a_2 & a_1 & a_0 & \ldots \\
\vdots & \vdots & \vdots & \ddots \end{array} \right] ,
\]
and its {\em symbol} $a(t)$ can be defined as a map from the unit circle to the
complex plane by
\[
a(t) := \sum_{j=- \infty}^{\infty} a_j t^j ,~~\mbox{ where}~
\sum_{j=- \infty}^{\infty} | a_j | < \infty .
\]
We will be interested in symbols $b(t)$ with only finitely many nonzero
Fourier coefficients:
\[
b(t) := \sum_{j=-r}^s b_j t^j ,~~t \in \mbox{ unit circle},
\]
and we will denote the corresponding infinite Toeplitz matrix as $T(b)$ and its
top $n$ by $n$ block as $T_n (b)$.
Assume that $b(t) \neq 0$ for all $t$ on the unit circle, and let $k$ be the 
{\em winding number} of $b$; i.e., the number of times that $b(t)$ surrounds
the origin in a counterclockwise direction as $t$ traverses the unit circle.
Then, according to the splitting theorem \cite{RochSilb}, as $n \rightarrow \infty$, 
the $|k|$ smallest singular values of $T_n (b)$ approach $0$ exponentially, while, 
for $n$ sufficiently large, the remaining singular values are bounded away from $0$ 
by a constant $d > 0$ that depends only on $b$.

\begin{figure}[ht]
\centerline{\epsfig{file=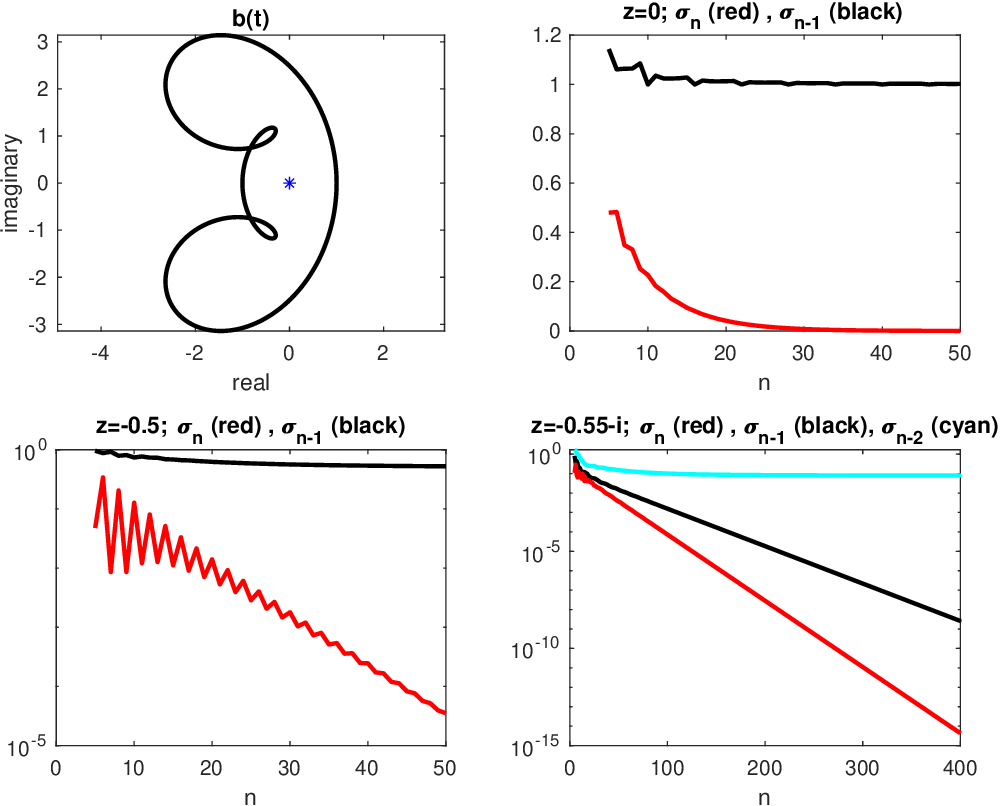,width=3in}}
\caption{Upper left plot shows image of the unit circle under 
$b(t) =- t^{-1} -1 + t + t^2 + t^3$.
The upper right and lower left plots show the smallest (red) and second smallest (black)
singular values of $T_n (b) - zI$ as a function of $n$ for $z=0$ and $z=-0.5$;
for these values, the winding number of $b(t) - z$ about the origin is $1$.
Such values lie in $S_{\epsilon} ( T_n (b) )$ for $n$ sufficiently large.
The lower right plot shows the smallest (red), second smallest (black), and
third smallest (cyan) singular values of $T_n (b) - zI$ as a function of $n$
for $z = -0.55 -i$.  In this case, the winding number of $b(t) - zI$ about the
origin is $2$, so the two smallest singular values go to $0$ as $n \rightarrow
\infty$, while the third smallest singular value remains bounded away from $0$.}
\label{fig:toeplitzexample}
\end{figure}

To determine $S_{\epsilon} ( T_n (b) )$ for $n$ sufficiently large, we must consider
the symbol $b(t) - z$ and look for values of $z$ for which the winding number of $b(t)-z$ 
about the origin is $\pm 1$.  The Grcar matrix, considered previously, is a Toeplitz
matrix whose symbol is $- t^{-1} + 1 + t + t^2 + t^3$.  By result 1 of 
Theorem \ref{thm:simple}, we can consider this matrix minus $2I$, so that the 
symbol is $b(t) = - t^{-1} -1 + t + t^2 + t^3$.  The image of the unit circle under $b$
is shown in the upper left plot of Figure \ref{fig:toeplitzexample}.  Clearly, the
absolute value of the winding number of this curve about the origin is $1$, and any
shift $z$ such that $b(t)-z$ still surrounds the origin with winding number $\pm 1$ will be in 
$S_{\epsilon} ( T_n (b) )$ for $n$ sufficiently large since the ratio of the smallest
to the second smallest singular value of $T_n (b) - zI$ will approach $0$ as 
$n \rightarrow \infty$.  This is illustrated in the upper right and lower left plots of 
Figure \ref{fig:toeplitzexample}, for $z = 0$ and $z =-0.5$.
These plots show the smallest (red) and second smallest (black) singular values
for $n=5$ to $50$.  The upper right plot is on a linear scale, while the lower
left plot uses a log scaling on the vertical axis to clearly show 
$\sigma_n ( T_n (b) - zI )$ approaching $0$ exponentially.  It seems clear from the plots that 
$\sigma_n ( T_n (b) - zI ) \rightarrow 0$ as $n \rightarrow \infty$ and that 
$\sigma_{n-1} ( T_n (b) - zI )$ approaches a nonzero value as $n \rightarrow \infty$.  
In the lower right plot, we chose a shift $z = -0.55-i$, so that the origin was
inside one of the small loops of the shifted symbol.  Here the winding number is $2$,
so the smallest and second smallest singular values of $T_n (b) - zI$ go to $0$
as $n \rightarrow \infty$, and it is the third smallest singular value that is 
bounded away from $0$.  This is illustrated in the lower right plot of 
Figure \ref{fig:toeplitzexample}.  In this case, the splitting theorem does
not guarantee that $z$ will lie in $S_{\epsilon} (A)$ for $n$ sufficiently
large, but it appears from the plot that this will be the case since 
$\sigma_n ( T_n (b) - zI )$ seems to be going to $0$ at a faster rate than
$\sigma_{n-1} ( T_n (b) - zI )$. 

\section{More General Analysis} \label{sec:genl}
Let $\lambda$ be a simple eigenvalue of $A$ and define $A_0 := A - \lambda I$.
Then $A_0$ has a zero eigenvalue that is also a zero singular value $\sigma_n ( A_0 )$.
The following theorem bounds the growth rate of $\sigma_n ( A_0 - zI )$ with $|z|$,
based on near-orthogonality properties that may be satisfied by the singular
vectors $u_n$ and $v_n$ corresponding to the zero singular value of $A_0$.

\begin{theorem}
\label{thm:genl}
Let $\lambda$ be a simple eigenvalue of $A$ and let $A_0 := A - \lambda I = U \Sigma V^{*}$
be a singular value decomposition of $A_0$, where $U := [ u_1 , \ldots , u_n ]$,
$V := [ v_1 , \ldots , v_n ]$, $\Sigma := \mbox{diag} ( \sigma_1 , \ldots , \sigma_{n-1} , 0 )$,
$\sigma_1 \geq \ldots \geq \sigma_{n-1} > 0$.  Let $A_0^{\dag}$ denote the Moore-Penrose
pseudoinverse of $A_0$:  
\begin{equation}
A_0^{\dag} := V_{n-1} \Sigma_{n-1}^{-1} U_{n-1}^H , \label{pseduoinverse}
\end{equation}
where $U_{n-1} := [ u_1 , \ldots , u_{n-1} ]$, $V_{n-1} := [ v_1 , \ldots , v_{n-1} ]$, and
$\Sigma_{n-1} = \mbox{diag} ( \sigma_1 , \ldots , \sigma_{n-1} )$.
For each $k=1,2, \ldots$, the smallest
singular value of $A_0 - zI$ is less than or equal to
\begin{equation}
\left\| \sum_{j=1}^k z^j ( u_n^H ( A_0^{\dag} )^{j-1} v_n )~ u_n + 
z^{k+1} ( A_0^{\dag} )^k~ v_n \right\| \leq
\label{sigmanbound2a}
\end{equation}
\begin{equation}
|z| \left( \sum_{j=1}^k  |z |^{j-1} | u_n^H ( A_0^{\dag} )^{j-1} v_n | + | z |^{k} / \sigma_{n-1}^k \right) .
\label{sigmanbound2}
\end{equation}
Assuming that the spectral radius of $z A_0^{\dag}$ is less than one, taking $k= \infty$
in (\ref{sigmanbound2a}) results in the following upper bound on the smallest
singular value of $A_0 - zI$:
\begin{equation}
|z|~ | u_n^H (I - z A_0^{\dag} )^{-1} v_n | . \label{sigmanbound2b}
\end{equation}
\end{theorem}

\begin{proof}
Define
\[
w := v_n + \sum_{j=1}^k z^j ( A_0^{\dag} )^j v_n .
\]
Then $\| w \| \geq 1$, since the range of $A_0^{\dag}$ is orthogonal to $v_n$.
Since $A_0 A_0^{\dag} = A_0 ( V_{n-1} \Sigma_{n-1}^{-1} U_{n-1}^H ) = U_{n-1} U_{n-1}^H$,
we can evaluate the product $( A_0 - zI )w$ term-by-term as:
\begin{eqnarray}
( A_0 - zI ) w & = & ( A_0 - zI ) v_n + \sum_{j=1}^k \left( z^j U_{n-1} U_{n-1}^H ( A_0^{\dag} )^{j-1} v_n -
z^{j+1} ( A_0^{\dag} )^j v_n \right) \nonumber \\
 & = & -z v_n - \sum_{j=2}^{k+1} z^j ( A_0^{\dag} )^{j-1} v_n + \sum_{j=1}^k z^j U_{n-1} U_{n-1}^H
( A_0^{\dag} )^{j-1} v_n \nonumber \\
 & = & - \sum_{j=1}^k z^j u_n u_n^H ( A_0^{\dag} )^{j-1} v_n - z^{k+1} ( A_0^{\dag } )^k v_n .
\label{A0mzIw}
\end{eqnarray}
The upper bound in (\ref{sigmanbound2a}) is the norm of $( A_0 - zI )w$, and that 
in (\ref{sigmanbound2})
is obtained by replacing each term in (\ref{A0mzIw}) by its norm and noting that
the norm of $( A_0^{\dag} )^k v_n$ is less than or equal to $1/ \sigma_{n-1}^k$.

Assuming that the spectral radius of $z A_0^{\dag}$ is less than one, the quantity
$z^{k+1} ( A_0^{\dag} )^k v_n$ in (\ref{A0mzIw}) will go to $0$ as $k \rightarrow \infty$,
and the summation in (\ref{A0mzIw}) becomes $-z ( u_n^H ( I - z A_0^{\dag} )^{-1} v_n )~ u_n$,
whose norm is given in (\ref{sigmanbound2b}).
\end{proof}

Taking $k=1$ in (\ref{sigmanbound2}), we obtain the bound
\[
\sigma_{n} ( A_0 -zI ) \leq | u_n^H v_n |~|z| + | z |^2 / \sigma_{n-1} .
\]
If $\lambda$ is an ill-conditioned eigenvalue, then we know that $| u_n^H v_n |$ is tiny
and if $\sigma_{n-1} \approx 1$, then this shows that $\sigma_n ( A_0 - zI )$ grows more
like $| z |^2$ than like $| z |$ for $| u_n^H v_n | << | z | << 1$.  
If $u_n$ is also nearly orthogonal to $A_0^{\dag} v_n$, then taking $k=2$
in (\ref{sigmanbound2}) suggests that the growth rate of $\sigma_n (A_0 -zI)$ with $|z|$ may be 
more like $| z |^3$, and the more powers $j$ for which 
$| u_n^H ( A_0^{\dag} )^j v_n |$ is small, 
the higher the power of $|z|$ describing the growth of $\sigma_n ( A_0 - zI )$,
for $|z| << 1$.  
The smallest bound in Theorem \ref{thm:genl} may be obtained for 
$k = \infty$, if $u_n$ is almost orthogonal to $(I - z A_0^{\dag} )^{-1} v_n$.

Note that the vector $w$ used in the proof of Theorem \ref{thm:Jordan}, namely,
$w = [ 1, z, z^2 , \ldots , z^{n-1} ]$, is of the same form as that used in the
proof of Theorem \ref{thm:genl}, since the singular value decomposition of a 
Jordan block can be written as $U \Sigma V^H$, where $U$ is the identity,
$\Sigma = \mbox{diag} ( 1, \ldots , 1, 0 )$, and
\[
V = \left[ \begin{array}{cccc}
0          & \ldots & 0 & 1 \\
1          & \ldots & 0 & 0 \\
  & \ddots &        & \vdots \\
  &        & 1 & 0 \end{array} \right] .
\]
Thus $v_n$ is the first unit vector, $u_n$ is the $n$th unit vector, and
the pseudoinverse is 
\[
A_0^{\dag} = \left[ \begin{array}{cccc}
0 & 0 & \ldots & 0 \\
1 & 0 & \ldots & 0 \\
  & \ddots & & \vdots \\
  &        & 1 & 0 \end{array} \right] .
\]
It follows that $( A_0^{\dag} )^j v_n$ is the $(j+1)$st unit vector, for 
$j=1, \ldots , n-1$ and $u_n$ is orthogonal to $( A_0^{\dag} )^j v_n$ for
$j=0, \ldots , n-2$, while $u_n^H ( A_0^{\dag} )^{n-1} v_n = 1$ and
$u_n^H ( I - z A_0^{\dag} )^{-1} v_n = z^{n-1}$.
In the next section, we will show how the pseudoinverse comes into the 
derivative of the singular vectors $v_n$ and $u_n$.

In all of the matrices associated with Figures \ref{fig:grcar50}--\ref{fig:sampling10}, 
$u_n$ is nearly orthogonal not just to $v_n$ but also to several powers 
$( A_0^{\dag} )^j v_n$.  Defining $A_0$ to be $( A - \lambda I )$, 
where $\lambda$ is an eigenvalue of $A$ with maximal condition
number, we computed $\sigma_{n-1}$, $| u_n^H v_n |$, $| u_n^H ( A_0^{\dag} )^j v_n |$, and 
$| u_n^H ( A_0^{\dag} )^j v_n | / \| ( A_0^{\dag} )^j v_n \|$ for $j=1, \ldots , 5$.  
For the Grcar matrix, $\sigma_{n-1} \approx 0.81$, 
$| u_n^H v_n | \approx 4.6e-8$, and $| u_n^H ( A_0^{\dag} )^j v_n |$, $j=1, \ldots , 5$, 
ranged from $1.0e-7$ to $5.7e-3$; after dividing by $\| ( A_0^{\dag} )^j v_n \|$ to measure
the level of orthogonality, the results ranged from $1.1e-7$ to $4.0e-3$.
For the transient\_demo matrix, $\sigma_{n-1} \approx 0.18$,
$| u_n^H v_n | \approx 1.4e-6$, and the values of 
$| u_n^H ( A_0^{\dag} )^j v_n |$ ranged from $2.7e-5$ to $2.5$, but there was still
a significant degree of orthogonality between the vectors, as the values of
$| u_n^H ( A_0^{\dag} )^j v_n | / \| ( A_0^{\dag} )^j v_n \|$ ranged from 
$5.0e-6$ to $4.3e-3$.  For the sampling matrix,
$\sigma_{n-1} \approx 6.6$, $| u_n^H v_n | \approx 7.6e-7$,
and the values of $| u_n^H ( A_0^{\dag} )^j v_n |$ ranged from $1.1e-6$ (for $j=2$)
down to $7.3e-8$ (for $j=5$),
but after normalizing, we see that, in addition to $u_n$ being nearly orthogonal
to $( A_0^{\dag} )^j v_n$, the norms of these vectors decreased with $j$;
the range of values for $| u_n^H ( A_0^{\dag} )^j v_n | / \| ( A_0^{\dag} )^j v_n \|$
was $3.0e-6$ (for $j=1$) to $1.3e-2$ (for $j=5$).

\begin{corollary}
\label{cor:bounds}
With the notation and assumptions of Theorem \ref{thm:genl}, 
let $\epsilon \in (0,1)$ be given.  The region $S_{\epsilon} (A)$ contains
the set of points $z \in \mathbb{C}$ such that $|z- \lambda | < \sigma_{n-1}$
and the expression in either (\ref{sigmanbound2a}), (\ref{sigmanbound2}), or
(\ref{sigmanbound2b}), with $z$ replaced by $(z- \lambda )$, is bounded above 
by $\epsilon ( \sigma_{n-1} - | z- \lambda | )$.
\end{corollary}
\begin{proof}
From Theorem \ref{thm:genl}, $\sigma_n ( A - zI) = 
\sigma_n ( A_0 - (z - \lambda ) I )$ is bounded above by the
quantities in (\ref{sigmanbound2a}), (\ref{sigmanbound2}), and (\ref{sigmanbound2b}),
when $z$ is replaced by $z- \lambda$ in these expressions.  
It follows from (\ref{sigmanm1bound}) that $\sigma_{n-1} ( A-zI) =
\sigma_{n-1} ( A_0 - (z - \lambda ) I)$ is bounded below by
$\sigma_{n-1} - |z - \lambda |$.  The result follows since
$\sigma_2 ( ( A-zI )^{-1} ) / \sigma_1 ( ( A-zI )^{-1} ) =
\sigma_{n-1} ( A-zI ) / \sigma_n (A-zI)$.
\end{proof}

\begin{figure}[ht]
\centerline{\epsfig{file=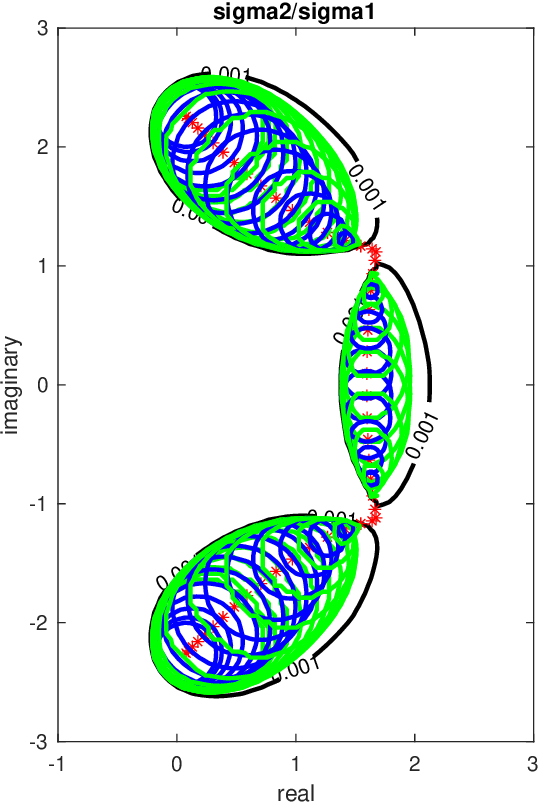,width=1.5in}}
\caption{Disks about each eigenvalue where inequality (\ref{ineq1}) is satisfied 
(outlined in blue) and regions about each eigenvalue where inequality (\ref{ineq2}) 
is satisfied (outlined in green).  $S_{10^{-3}} (A)$ is outlined in black.
Matrix is the Grcar matrix of size $n=50$.}
\label{fig:grcar50bounds}
\end{figure}

\begin{figure}[ht]
\centerline{\epsfig{file=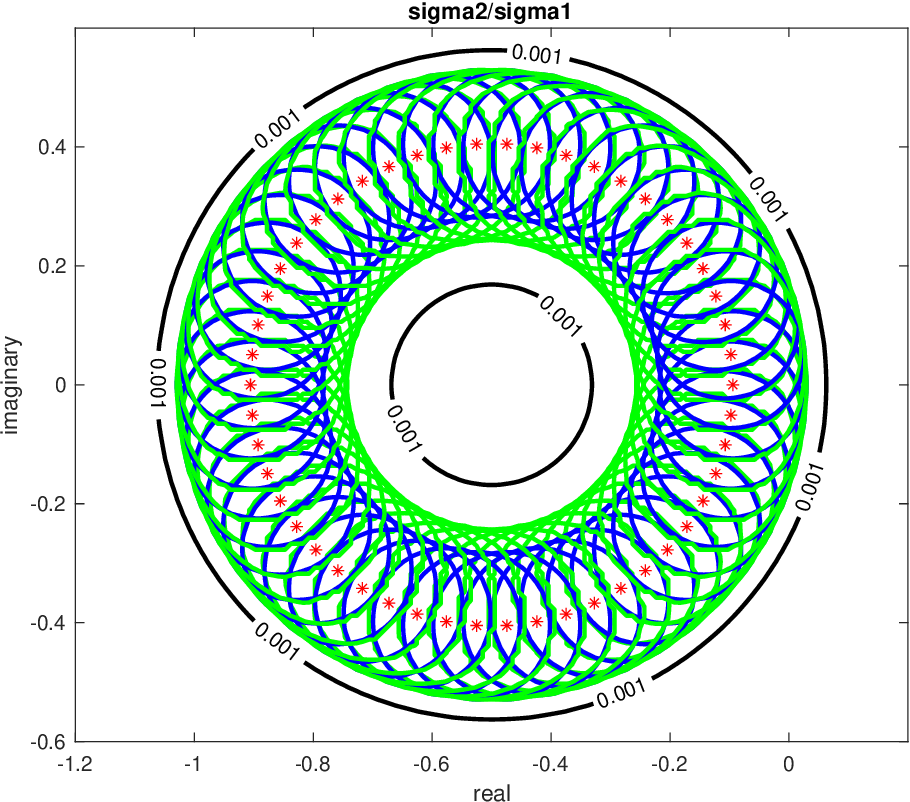,width=2in}}
\caption{Disks about each eigenvalue where inequality (\ref{ineq1}) is satisfied
(outlined in blue) and regions about each eigenvalue where inequality (\ref{ineq2})
is satisfied (outlined in green).  $S_{10^{-3}} (A)$ is outlined in black. 
Matrix is the transient\_demo matrix of size $n=50$.}
\label{fig:transient50bounds}
\end{figure}

\begin{figure}[ht]
\centerline{\epsfig{file=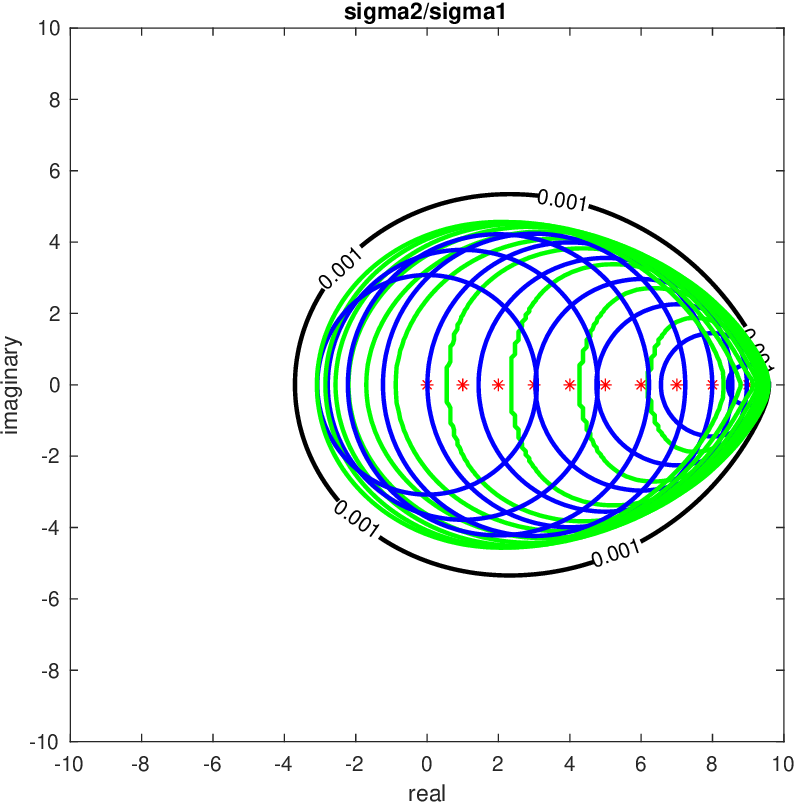,width=2in}}
\caption{Disks about each eigenvalue where inequality (\ref{ineq1}) is satisfied
(outlined in blue) and regions about each eigenvalue where inequality (\ref{ineq2})
is satisfied (outlined in green).  $S_{10^{-3}} (A)$ is outlined in black. 
Matrix is the sampling matrix of size $n=10$.}
\label{fig:sampling10bounds}
\end{figure}

Figures \ref{fig:grcar50bounds}--\ref{fig:sampling10bounds} show the 
regions about each eigenvalue where 
\begin{equation}
\min_{k=1, \ldots , 100} | z- \lambda | \left( \sum_{j=1}^k | z - \lambda |^{j-1} 
| u_n^H ( A_0^{\dag} )^{j-1} v_n | + | z - \lambda |^k / \sigma_{n-1}^k \right) \leq
10^{-3} ( \sigma_{n-1} - | z - \lambda | ) , \label{ineq1}
\end{equation}
and where
\begin{equation}
| z - \lambda |~ | u_n^H ( I - (z- \lambda ) ( A_0^{\dag} )^{-1} v_n | \leq 
10^{-3} ( \sigma_{n-1} - | z - \lambda | ) , \label{ineq2}
\end{equation}
for the Grcar matrix, the transient\_demo matrix and the sampling matrix.
Note that inequality (\ref{ineq1})  describes disks about each eigenvalue (outlined in blue
in the figures) since it depends only on $|z - \lambda |$, while inequality (\ref{ineq2}) 
describes more general regions (outlined in green). 
The disks outlined in blue cover a good portion of $S_{10^{-3}} (A)$, and the 
regions outlined in green cover still more of $S_{10^{-3}} (A)$.

Theorem \ref{thm:genl} also leads to a new result about pseudospectra:

\begin{corollary}
\label{cor:pseudo}
With the notation and assumptions of Theorem \ref{thm:genl}, 
let $\epsilon > 0$ be given.  The $\epsilon$-pseudospectrum contains
the set of points $z \in \mathbb{C}$ such that
the expression in (\ref{sigmanbound2a}), (\ref{sigmanbound2}), or 
(\ref{sigmanbound2b}), with $z$ replaced by $z- \lambda$, is bounded above
by $\epsilon$.
\end{corollary}
\begin{proof}
From Theorem \ref{thm:genl}, $\sigma_n ( A - zI) = 
\sigma_n ( A_0 - (z - \lambda ) I )$ is bounded above by the
expressions in (\ref{sigmanbound2a}), (\ref{sigmanbound2}), and (\ref{sigmanbound2b}),
when $z$ is replaced by $z - \lambda$ in these expressions.
The result follows since $1/ \sigma_1 ( ( A-zI )^{-1} ) = \sigma_n (A-zI)$.
\end{proof} 

\begin{figure}[ht]
\centerline{\epsfig{file=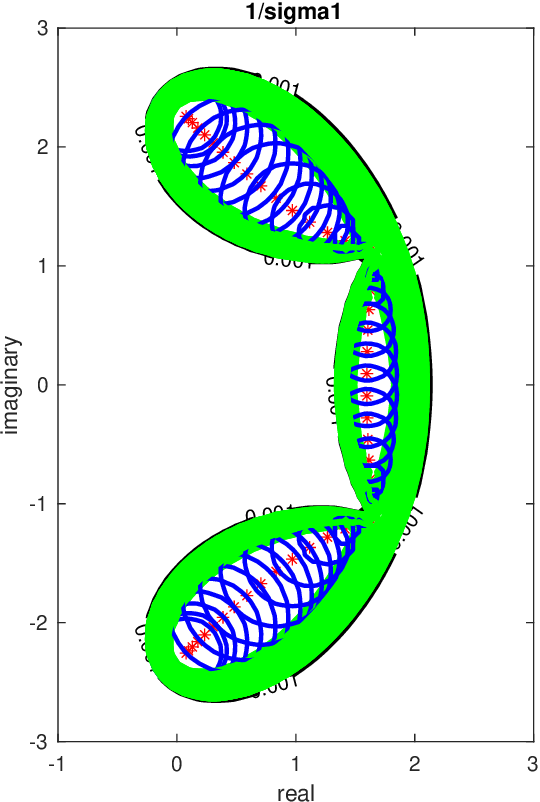,width=2in}}
\caption{Disks about each eigenvalue derived from inequality (\ref{sigmanbound2})
(outlined in blue) and regions about each eigenvalue derived from inequality
(\ref{sigmanbound2b}) (outlined in green) that lie in the $10^{-3}$-pseudospectrum
(outlined in black).  Matrix is the Grcar matrix of size $n=50$.}
\label{fig:grcar50pseudodisks}
\end{figure}

The regions defined in Corollary \ref{cor:pseudo} come even closer to filling
the $\epsilon$-pseudospectrum of $A$.
We illustrate this just for the Grcar matrix in Figure \ref{fig:grcar50pseudodisks}.
Interestingly, each of the regions outlined in green, where
$| z - \lambda |~ | u_n^H ( I - (z- \lambda ) ( A_0^{\dag} )^{-1} v_n | \leq 10^{-3}$,
contains a large portion of $\Lambda_{10^{-3}} (A)$, and not just a region
around the eigenvalue $\lambda$.

\section{Differentiation of Singular Values and Singular Vectors} \label{sec:diffeq}

In the previous section, we described regions where $\sigma_2 (( A-zI )^{-1} ) /
\sigma_1 (( A-zI )^{-1} ) < \epsilon$.  We now wish to consider regions
where $| u_n (z )^H v_n (z) | < \epsilon$, where $u_n (z)$ and $v_n (z)$
are right and left singular vectors associated with the largest singular
value of $(A-zI )^{-1}$, that is, with the smallest singular value of $A-zI$.
It was noted in Section \ref{sec:Numerical} that these regions appear to be 
slightly smaller.

In general, the SVD is {\em not} continuously differentiable unless singular
values are allowed to become negative \cite{deMoor1989,Bunse1991,Wright1992}.
However, one can differentiate in one direction, provided there are no multiple
singular values (and, in some cases, even if there are multiple singular values)
\cite{ONeil}.  Therefore, let us fix $\theta \in [0, 2 \pi )$, set
$A(r) := A_0 - r e^{i \theta} I$, $r \geq 0$, where $A_0$ has a simple eigenvalue
$0$ that is also a singular value, and write the SVD of $A(r)$ as
$A(r) = U(r) \Sigma (r) V(r )^H$.  Then
\begin{equation}
A' (r) = - e^{i \theta} I = U' (r) \Sigma(r) V(r )^H + U(r) \Sigma' (r) V(r )^H +
U(r) \Sigma (r) ( V' (r) )^H . \label{SVDderiv}
\end{equation}
We also require $U(r )^H U(r) = V(r )^H V(r) = I$, and differentiating these
equations leads to
\begin{equation}
( U' (r) )^H U(r) + U(r )^H U' (r) = 
( V' (r) )^H V(r) + V(r )^H V' (r) = 0 . \label{UVderiv}
\end{equation}
Multiplying (\ref{SVDderiv}) on the left by $U(r )^H$ and on the right by $V(r)$
leads to the relation
\begin{equation}
- e^{i \theta} U(r )^H V(r) = U(r )^H U' (r) \Sigma (r) + \Sigma' (r) +
\Sigma (r) ( V' (r) )^H V(r) . \label{relation}
\end{equation}
Equating real parts of the diagonal entries in (\ref{relation}) and noting that 
equation (\ref{UVderiv}) says that $U(r )^H U'(r)$ and $( V' (r) )^H V(r)$ are 
skew Hermitian so that the real parts of their diagonal entries are zero, we find
\begin{equation}
- \mbox{Re} ( e^{i \theta} \mbox{diag} ( U(r )^H V(r) )) = \mbox{diag} ( \Sigma' (r) ) .
\label{Sigmap}
\end{equation}
Thus the derivative of $\sigma_n (r)$ is
\begin{equation}
\sigma_n' (r) = - \mbox{Re} ( e^{i \theta} u_n (r )^H v_n (r) ) ,
\label{sigmanp}
\end{equation}
showing that a small inner product $| u_n (r )^H v_n (r) |$ between
left and right singular vectors leads to slow change in the singular value
$\sigma_n (r)$.  Conversely, if one can establish that $\sigma_n (r)$
changes slowly with $r$ in all directions, then this would imply that
$| u_n (r )^H v_n (r) |$ must be small; unfortunately, in previous sections
we have established only upper bounds on $\sigma_n (r)$, leaving open the
possibility that it might alternate quickly between $0$ and these upper bounds.

Defining $Z(r) := U(r )^H U' (r)$ and $W(r) := V(r )^H V' (r)$ and
equating imaginary parts of the diagonal entries in (\ref{relation}) leads
to $n$ equations in $2n$ unknowns:
\begin{equation}
- i~ \mbox{Im} ( e^{i \theta} \mbox{diag} ( U(r )^H V(r) )) = 
\mbox{diag} ( Z(r) ) \Sigma (r) + \Sigma (r) 
\mbox{diag} ( W(r )^H )  . \label{imagdiag}
\end{equation}
Any solution to these equations, say, taking $\mbox{diag} ( Z(r) ) = 
\mbox{diag} ( W(r )^H )$ should lead to the same result that $A(r) =
U(r) \Sigma (r) V(r )^H$, where $U(r )^H U(r) = V(r )^H V(r) = I$.
Note, however, that if $\sigma_n (0) = 0$, then these equations have no 
solution unless $\mbox{Im} ( e^{i \theta} u_n (0 )^H v_n (0) ) = 0$.
In this case, we can replace any computed $u_n (0)$ by $e^{i \varphi} u_n (0)$,
where $\varphi = \theta + \mbox{arg} ( u_n (0 )^H v_n (0) ) + \pi$ so that
$\sigma _n ' (0) = - \mbox{Re} ( e^{i \theta} u_n (0 )^H v_n (0) )$ will
be nonnegative.

Equating the $n^2 - n$ off-diagonal entries in (\ref{relation}) and using
the fact that $Z(r)$ and $W(r)$ are skew Hermitian, we can solve for $Z(r)$
and $W(r)$.  Define $Q(r) := U(r )^H A' (r) V(r) = - e^{i \theta} U(r )^H V(r)$. 
Dropping the argument $r$ for convenience, and looking at the $(j,k)$ and $(k,j)$ 
entries of (\ref{relation}), we find
\begin{eqnarray}
\sigma_k z_{jk} + \sigma_j \bar{w}_{kj} & = & Q_{jk} , \label{firsteq} \\
\sigma_j z_{kj} + \sigma_k \bar{w}_{jk} & = & Q_{kj} . \nonumber
\end{eqnarray}
Substituting $z_{kj} = - \bar{z}_{jk}$ and $\bar{w}_{jk} =
- w_{kj}$, the second equation becomes
\[
\sigma_j \bar{z}_{jk} + \sigma_k w_{kj} = - Q_{kj} ,
\]
and taking complex conjugates
\begin{equation}
\sigma_j z_{jk} + \sigma_k \bar{w}_{kj} = - \overline{Q_{kj}} . \label{secondeq}
\end{equation}
Assuming that $\sigma_j \neq \sigma_k$,
equations (\ref{firsteq}) and (\ref{secondeq}) now enable us to solve
for $z_{jk}$ and $\bar{w}_{kj}$, say, for $j < k$ and we use the
fact that these matrices are skew-Hermitian to set the rest of
their off-diagonal entries.  To see what to do in case $\sigma_j = \sigma_k$,
see \cite{Wright1992}.
Finally, after solving for $Z(r)$ and $W(r)$, we get $U' (r) = U(r) Z(r)$
and $V' (r) = V(r) W(r)$.

We are interested in the quantity $u_n (r )^H v_n (r)$,
whose derivative is $( u_n' (r) )^H v_n (r) +
u_n (r )^H v_n' (r)$.  Since $u_n' (r) = U(r) z_n (r)$, where $z_n (r)$
is the $n$th column of $Z(r)$ and $v_n' (r) = V(r) w_n (r)$, where 
$w_n (r)$ is the $n$th column of $W(r)$, we can use formulas (\ref{firsteq})
and (\ref{secondeq}) with $k=n$ and $j < n$, along with the formulas
$Q_{jn} = - e^{i \theta} u_j^H v_n$ and $Q_{nj} = - e^{i \theta} u_n^H v_j$,
to solve for the relevant quantities.  Again dropping the argument $r$
for convenience, we can write
\begin{eqnarray}
z_{jn} & = & - ( e^{i \theta} \sigma_n u_j^{H} v_n + e^{-i \theta} \sigma_j v_j^{H} u_n ) /
( \sigma _n^2 - \sigma_j^2 ) , \label{zjn} \\
\bar{w}_{nj} & = & ( e^{-i \theta} \sigma_n v_j^{H} u_n + e^{i \theta} \sigma_j u_j^{H} v_n ) /
( \sigma_n^2 - \sigma_j^2 ) .  \label{wbarnj}
\end{eqnarray}
Since $W$ is skew Hermitian,
\begin{equation}
w_{jn} = - ( e^{-i \theta} \sigma_n v_j^{H} u_n + e^{i \theta} \sigma_j u_j^{H} v_n ) /
( \sigma_n^2 - \sigma_j^2 ) . \label{wjn}
\end{equation}
It follows that
\begin{eqnarray}
v_n^{H} u_n' = v_n^{H} U z_n & = & v_n^{H} \left[ \sum_{j=1}^{n-1}
\frac{e^{i \theta} \sigma_n ( u_j^{H} v_n ) + e^{-i \theta} \sigma_j ( v_j^{H} u_n )}
{\sigma_j^2 - \sigma_n^2} u_j  + z_{nn} u_n \right] \nonumber \\
 & = & e^{i \theta} \sigma_n \sum_{j=1}^{n-1} \frac{| u_j^{H} v_n |^2}
{\sigma_j^2 - \sigma_n^2} +
e^{-i \theta} \sum_{j=1}^{n-1} \sigma_j \frac{( v_j^{H} u_n ) ( v_n^{H} u_j )}
{\sigma_j^2 - \sigma_n^2} + z_{nn} v_n^{H} u_n .
\label{vnstarunp}
\end{eqnarray}
When $r=0$, the first term in (\ref{vnstarunp}) is $0$ since $\sigma_n (0) = 0$.
The third term involves $z_{nn} (0)$, but based on the argument after (\ref{imagdiag}), 
we can take $z_{nn} (0) = 0$.  Thus we are left with only the second term in
(\ref{vnstarunp}), which can be written in the form
\[
e^{-i \theta} v_n^H U_{n-1} \Sigma_{n-1}^{-1} V_{n-1}^H u_n =
e^{-i \theta} v_n^H ( A_0^{\dag} )^H u_n ,
\]
where again the pseudoinverse $A_0^{\dag}$ is involved.  Thus
\begin{equation}
( u_n' (0) )^H v_n (0) = e^{i \theta} u_n (0 )^H A_0^{\dag} v_n (0) . \label{unpstarvn0}
\end{equation}   

Similarly, we can write
\begin{eqnarray}
u_n^{H} v_n' = u_{n}^{H} V w_n & = & u_n^{H} \left[ \sum_{j=1}^{n-1}
\frac{e^{-i \theta} \sigma_n ( v_j^{H} u_n ) + e^{i \theta} \sigma_j ( u_j^{H} v_n )}
{\sigma_j^2 - \sigma_n^2} v_j + w_{nn} v_n \right] \nonumber \\
 & = & e^{-i \theta} \sigma_n \sum_{j=1}^{n-1} \frac{| v_j^{H} u_n |^2}
{\sigma_j^2 - \sigma_n^2} + e^{i \theta} \sum_{j=1}^{n-1} \sigma_j
\frac{( u_j^{H} v_n )( u_n^{H} v_j )}{\sigma_j^2 - \sigma_n^2} + w_{nn} u_n^{H} v_n .
\label{unstarvnp}
\end{eqnarray}
Again for $r=0$, the first term in (\ref{unstarvnp}) is $0$ since $\sigma_n (0) = 0$,
and the third term is $0$ since $w_{nn} (0)$ can be taken to be $0$.  Again
we are left with the second term in (\ref{unstarvnp}), which can be written as
\[
e^{i \theta} u_n^H V_{n-1} \Sigma_{n-1}^{-1} U_{n-1}^H v_n = 
e^{i \theta} u_n^H A_0^{\dag} v_n .
\]
Thus
\begin{equation} 
u_n (0 )^H v_n' (0) = e^{i \theta} u_n (0 )^H A_0^{\dag} v_n (0) . \label{unstarvnp0}
\end{equation}

Combining (\ref{unpstarvn0}) and (\ref{unstarvnp0}), we see that
\begin{equation}
\left. \frac{d}{dr} ( u_n (r )^H v_n (r) ) \right|_{r=0} =
2 e^{i \theta} u_n (0 )^H A_0^{\dag} v_n (0) . \label{innerprodderiv}
\end{equation}
Again the quantity $u_n (0 )^H A_0^{\dag} v_n (0)$ is related to the 
rate of change of $ u_n (r )^H v_n (r )$ with $r$, as it was related to the
rate of change of $\sigma_n ( A_0 - zI)$  in Theorem \ref{thm:genl}.  If
$| u_n (0 )^H A_0^{\dag} v_n (0) |$ is tiny then we expect
$| u_n (r )^H v_n (r) |$ to grow (or decay) more like $r^2$ than like $r$,
for small $r$.  Unfortunately, higher derivatives of $u_n (r )^H v_n (r)$
are significantly more complicated, and we do not have a result analogous
to that in Corollary \ref{cor:bounds} describing a region where 
$| u_n (r )^H v_n (r) | < \epsilon$.  This remains an open problem.

\section{Some Possible Applications and Remaining Open Problems} \label{sec:open}
A number of applications, most notably in fluid mechanics, involve working
with large resolvent matrices \cite{application1,application2,application3},
and it has been observed, based on physical grounds, that it may be sufficient
to work with a rank one or low rank approximation to the resolvent.
(For some plots of singular values of such resolvents, see 
\cite{resplot1,resplot2}.)
When a large resolvent matrix can be replaced by a rank one or low rank matrix, 
this results in a huge savings in storage and compute time.  In this paper, we
have used a different rank one approximation to $(A-zI )^{-1}$ for each $z$,
but if the same low rank approximation can be used for many different values
of $z$, then further savings is possible.

In \cite{GreenWell}, it was observed that two upper bounds on the norm
of an analytic function $f(A)$ turned out to be approximately the same because
the numerical range of the resolvent was close to a disk about the origin;
this led to the realization that the resolvent resembled a rank one matrix
$\sigma_1 (z) u_1 (z) v_1 (z )^H$, where $u_1 (z)$ and $v_1 (z)$ are almost 
orthogonal to each other.  With more knowledge of the behavior of the
singular vectors $u_1 (z)$ and $v_1 (z)$, one might be able to say more
about the sharpness, or lack thereof, of these bounds on $\| f(A) \|$.

It is surprising that, to the best of our knowledge, the connection
between pseudospectra and regions where the resolvent is close to a 
rank one matrix has not been explained before.  The proof just rests
on the fact that if the second smallest singular value of the matrices
$A_{0,j}$ in Theorem \ref{thm:relation} is on the order of $1$, then from
(\ref{sigmanm1bound}), the second smallest singular value of $A_{0,j} - zI$
is between $1 - |z|$ and $1 + |z|$.

In Theorem \ref{thm:genl}, we do not know how to relate the 
properties that $u_n$ be almost orthogonal to $( A_0^{\dag} )^j v_n$ 
to other standard matrix properties.  The property that $u_n$ be almost
orthogonal to $v_n$ just corresponds to $\lambda$ being ill-conditioned
since the condition number of $\lambda$ is $1/ | u_n^H v_n |$.  We also
do not know for what class of matrices these additional near orthogonality
properties hold.  They hold for the test problems presented in this 
paper and for many others that we have tried, but identifying this class
of matrices remains an open problem.

Finally, we still have not delineated the region where 
$| u_1 (z )^H v_1 (z) | < \epsilon$, although in the plots, it appears to
be most of $S_{\epsilon} (A)$.  Identifying this region remains an open problem.
 
\vspace{.1in}
\noindent
{\bf Acknowledgments.} The first author thanks Mark Embree for helpful discussions,
especially the suggestion to look at Toeplitz matrices.  The authors thank
the University of Washington Applied Math Department for hosting the sabbatical
visit of the second and third authors, during which time this work began.

\end{document}